\documentclass[12pt]{article}
\usepackage{latexsym}
\usepackage{amssymb}
\usepackage[dvips]{graphicx}
\usepackage{epsfig}
\usepackage{rotating}
\usepackage{amsfonts}
\usepackage{amsmath}

\usepackage{hyperref}
\usepackage{pstricks}

\newtheorem{theorem}{Theorem}[section]

\newtheorem{definition}{Definition}[section]

\newcommand{\beqa}{\begin{eqnarray}}
\newcommand{\eeqa}{\end{eqnarray}}

\newcommand{\bea}{\begin{eqnarray}}
\newcommand{\eea}{\end{eqnarray}}

\newcommand{\be}{\begin{equation}}
\newcommand{\ee}{\end{equation}}
\newcommand{\beq}{\begin{eqnarray}}
\newcommand{\eeq}{\end{eqnarray}}

\def\lsim{\
  \lower-2.0pt\vbox{\hbox{\rlap{$<$}\lower5.5pt\vbox{\hbox{$\sim$}}}}\ }
\def\gsim{\
  \lower-2.0pt\vbox{\hbox{\rlap{$>$}\lower5.5pt\vbox{\hbox{$\sim$}}}}\ }

\begin{document}

\begin{titlepage}
\begin{flushright}
\end{flushright}

\vspace{20pt}

\begin{center}

{\Large\bf 
Generalization of the Bollob\'as-Riordan polynomial for tensor graphs}
\vspace{20pt}

 Adrian Tanasa

\vspace{15pt}

\vspace{10pt}$^{}${\sl
LIPN, Institut Galil\'ee, CNRS UMR 7030,\\
Univ. Paris Nord, 99 av. Cl\'ement, 93430 Villetaneuse, France, UE}\\

\vspace{10pt}
$^{}${\sl
Departamentul de Fizic\u a Teoretic\u a,\\ 
Institutul de Fizic\u a \c si Inginerie Nuclear\u a Horia Hulubei,\\
P. O. Box MG-6, 077125 M\u agurele, Rom\^ania, UE}\\

\vspace{20pt}
E-mail:  
 {\em adrian.tanasa@ens-lyon.org}

\vspace{10pt}

\begin{abstract}
\noindent
Tensor models are used nowadays for implementing a fundamental theory of quantum gravity. We define here a polynomial $\mathcal T$ encoding the supplementary topological information. This polynomial is a natural generalization of the Bollob\'as-Riordan polynomial (used to characterize matrix graphs) and is different of the Gur\u au polynomial (R. Gur\u au, ``Topological Graph Polynomials in Colored Group Field Theory'', Annales Henri Poincare {\bf 11}, 565-584 (2010)), defined for a particular class of tensor graphs, the colorable ones. The polynomial  $\mathcal T$ is defined for both colorable and non-colorable  graphs and it is proved to satisfy the deletion/contraction relation. A non-trivial example of a non-colorable graphs is analyzed.
\end{abstract}

\today

\end{center}


\end{titlepage}


\setcounter{footnote}{0}

\section{Introduction and motivation}
\label{Intro}
\renewcommand{\theequation}{\thesection.\arabic{equation}}
\setcounter{equation}{0}

Tensor models are nowadays one of the candidates for a fundamental theory of quantum gravity. The theoretical physics framework is a quantum field theoretical (QFT) one, namely group field theory (GFT) (the interested reader can turn to the review papers \cite{gft}). Roughly speaking, the main idea is that if matrix models can be seen as dual to triangulations of two-dimensional surfaces, this can then be extended - using tensor models instead of matrix models - to three-dimensional and finally to four-dimensional spaces. It is worth emphasizing, that unlike string theories, these models do not require extra-dimensions, so that one can work in a space with dimension four, the dimension of space-time at our energy scale.

Let us also mention that an important regain of interest for GFTs is to be noticed lately in the mathematical and theoretical physics communities (see for example \cite{gft-regain} and references within).

In this paper we deal with the three-dimensional case, which is known to correspond to the topological version of three-dimensional quantum gravity (the GFT formulation of the Ponzano-Regge model).

The topological richness of three-dimensional tensor graphs is far from being completely tamed. Within this framework, an important notion is the one of {\it bubbles}, which are 
 a natural generalization of the notion of faces appearing in the case of ribbon graphs (also known as maps, see for example \cite{maps}).
The bubbles can thus be seen as ribbon graphs and one can simply calculate their genera (as Riemann surfaces). An important result is that, if all the bubbles are planar ({\it i. e.} they have vanishing genera), then the respective graph corresponds to a {\it manifold.} 

An useful tool for understanding graphs is, from a combinatorial point of view, the definition of an universal graph polynomial. Thus, the Tutte polynomial \cite{tutte} is known to encode the topological information of some given graph. For ribbon graphs (also known as maps), the Tutte polynomial generalizes in a very elegant way - the Bollob\'as-Riordan polynomial \cite{br}. Both these polynomial obey a crucial relation, the deletion/contraction one.

When uplifting to three-dimensional tensor graphs, a first attempt to define such a polynomial generalizing the Tutte and Bollob\'as-Riordan ones is the Gur\u au polynomial \cite{gp}. This polynomial introduces a set of supplementary variables to keep track of the supplementary topological information. Furthermore, in \cite{gp} a modified notion of graph contraction was introduced. The Gur\u au polynomial satisfies then the deletion/contraction relation. Nevertheless, let us also recall that the Gurau polynomial was defined in \cite{gp} only for a particular class of tensor graphs, the colorable graphs.

In this paper we propose a different polynomial, called $\mathcal T$ which generalizes in a natural way the Tutte and Bollob\'as-Riordan ones. Thus, we introduce one supplementary variable which keeps track of the sum of genera of the bubbles of the graph. This definition is  particularly sensitive to the manifold-like character of the respective graph.  Let us also emphasize that the polynomial  $\mathcal T$ is fundamentally different of the Gur\u au one. Finally, let us also stress on the fact that the deletion/contraction relation for tensor graphs is proved to hold also for the  $\mathcal T$ polynomial.

\bigskip

The paper is structured as follows. In the next section we recall the definitions and main properties of the Tutte and Bollob\'as-Riordan polynomials. The third section introduces GFTs, emphasizing some of the non-trivial issues of its diagrammatic, like for example the way of defining bubbles. The notion of colorability of graphs is also recalled. The following section presents the Gur\u au polynomial.  The fifth section defines the polynomial  $\mathcal T$ and proves that the deletion/contraction relation holds for any tensor graph, colorable or not. The last section analyzes a non-trivial example, namely some non-colorable graph.

\section{The Tutte and the Bollob\'as-Riordan polynomials}
\label{TBR}
\renewcommand{\theequation}{\thesection.\arabic{equation}}
\setcounter{equation}{0}

We briefly recall here the definitions of the Tutte and Bollob\'as-Riordan polynomials as well as the deletion/contraction relation for both of them. The interested reader can turn to \cite{io-BR} for details.

\medskip

A graph $G$ is defined as a set of vertices $V$ and of edges $E$ together with 
an incidence relation between them. The number of vertices and edges in a graph
will be noted also $V$ and $E$ for simplicity, since our context always prevents any confusion. 

\begin{definition}
\noindent
\begin{enumerate}
\item A graph edge which is neither a bridge nor a self-loop is called {\bf regular}.
\item The rank of a subgraph $A$ ($A\subset G$) is defined as
\beqa
r(A)=V-k(A),
\eeqa
where $k(A)$ is the number of connected components of the subgraph $A$.
\item The nullity (or cyclomatic number) of a subgraph $A$ is defined as
\beqa
n(A)=\vert A\vert - r(A).
\eeqa
\end{enumerate}
\end{definition}

\begin{definition}
 If $G$ is a graph, then its Tutte polynomial of
$T_G(x,y)$ 
is defined as
\be
T_G (x,y)=\sum_{A\subset G}      (x-1)^{r(G)-r(A)} (y-1)^{n(A)}.
\ee
\end{definition}

\begin{figure}
\begin{center}
\includegraphics[scale=0.3,angle=0]{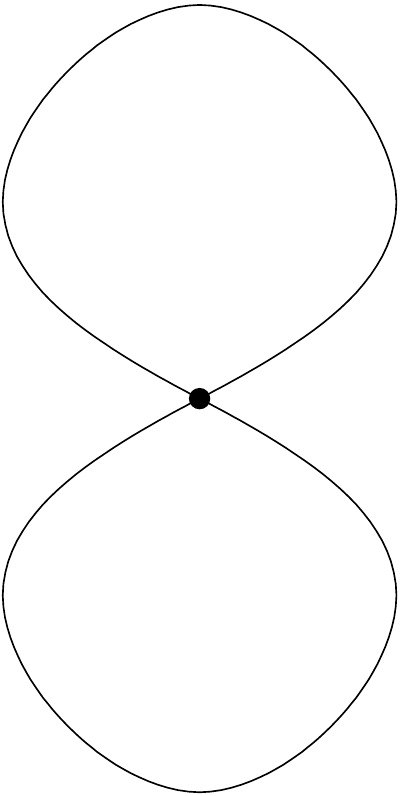}
\caption{The graph $G/e$.}
\label{fig:G/e}
\end{center}
\end{figure}

Let us now define two natural operations associated to some arbitrary edge $e$ of a graph $G$:
\begin{enumerate}
\item
The first of them is the {\it deletion}; it leads to a graph $G-e$. Thus, for the graph of Fig. \ref{fig:G}, deleting the edge $e$ leads to the graph $G-e$ represented in Fig.  \ref{fig:G-e}.

\begin{figure}
\begin{center}
\includegraphics[scale=0.3,angle=0]{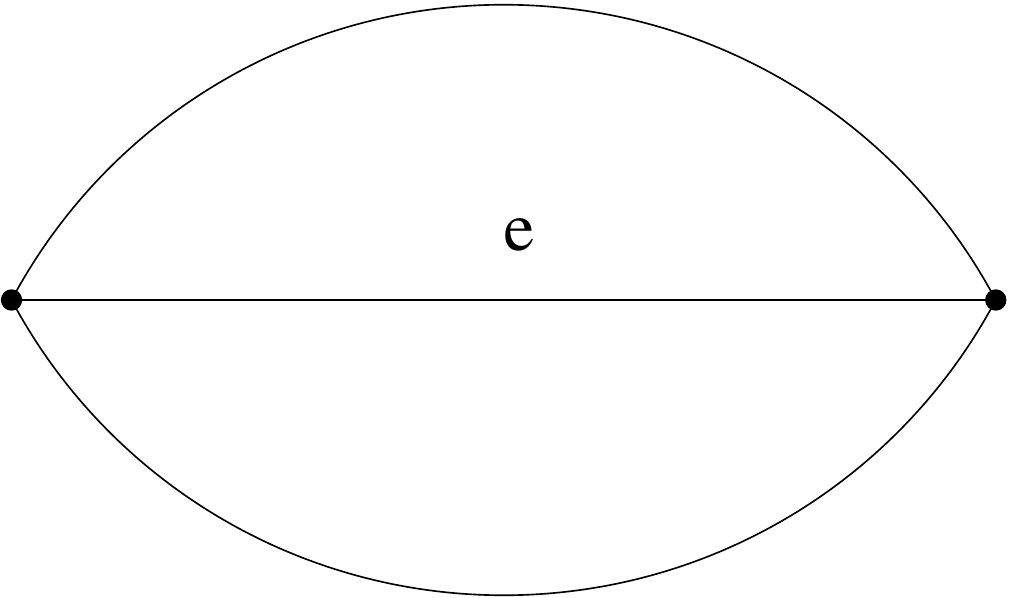}
\caption{An example of a graph $G$, with a chosen edge $e$.}
\label{fig:G}
\end{center}
\end{figure}

\begin{figure}
\begin{center}
\includegraphics[scale=0.3,angle=0]{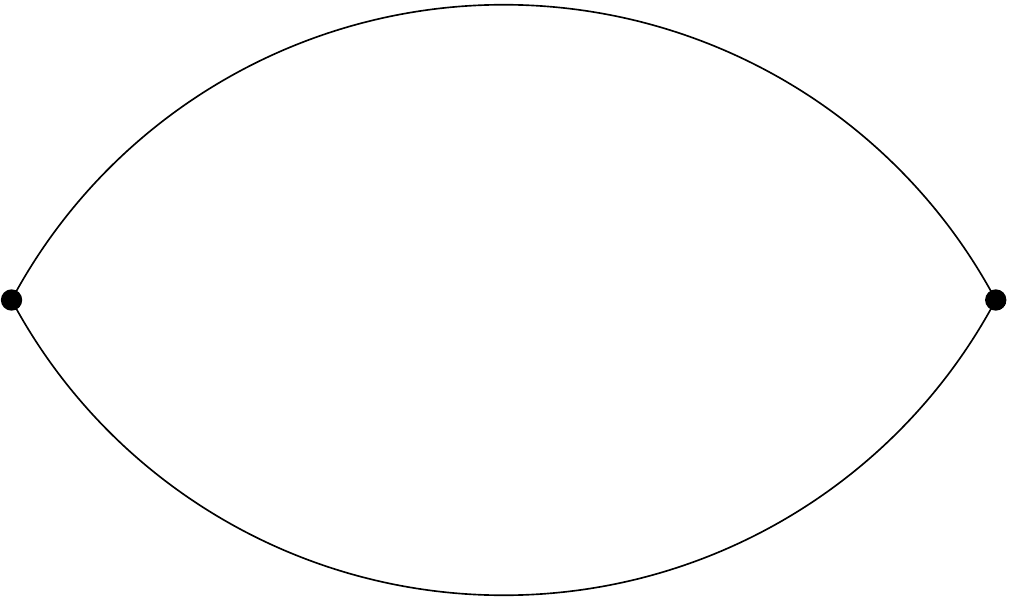}
\caption{The graph $G-e$.}
\label{fig:G-e}
\end{center}
\end{figure}

\item The second operation is the {\it contraction}. It deletes the respective edge $e$ and it contracts the two vertices that the $e$ was hooked to. The resulting graph is denoted $G/e$. For the example of Fig.  \ref{fig:G}, when one contracts the edge $e$ the graph $G/e$ is represented in Fig. \ref{fig:G/e}.   
\end{enumerate}

Using these definitions, we can now state the deletion/contraction property of the Tutte polynomial:

\begin{theorem}
\label{defdelcontr1}
If $G$ is a graph, and $e$ is a regular edge, then
\begin{equation}
T_G (x,y)=T_{G/e} (x,y)+T_{G-e} (x,y).
\end{equation}
\end{theorem}

\begin{figure}
\begin{center}
\includegraphics[scale=0.8,angle=-90]{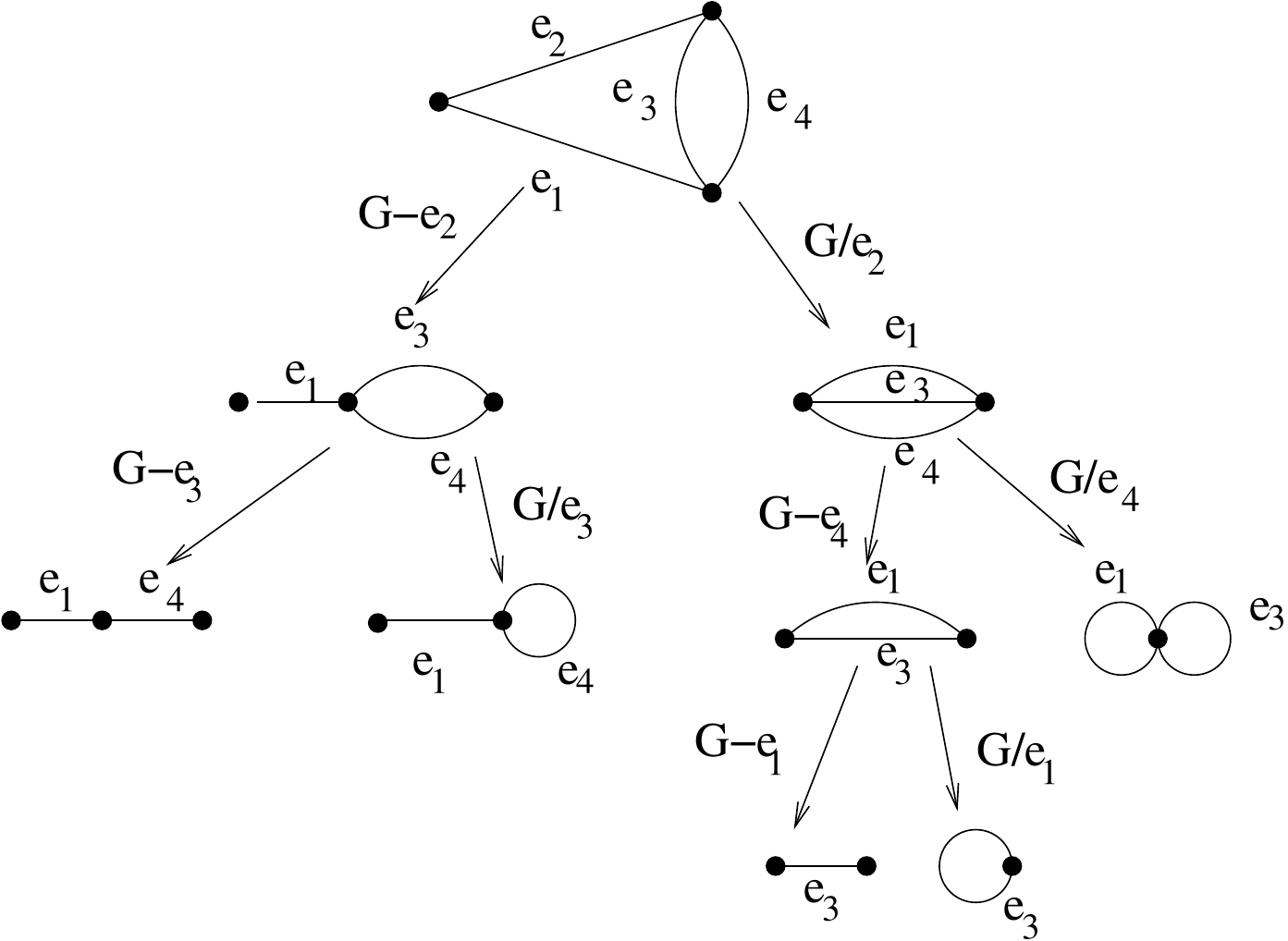}
\caption{The deletion/contraction of the regular edges of some graph $G$.}
\label{fig:cond}
\end{center}
\end{figure}

As an illustration for this kind of manipulation, an explicit example is the one of Fig. \ref{fig:cond}.

\medskip

A multivariate version of the Tutte polynomial exists in the literature \cite{sokal} (see also \cite{j0}). The main idea is that instead of having a single variable, $y$, for the number of edges, one introduces a set of variables $\beta_1,\ldots, \beta_E$, one for each edge:

\be
Z_G (q,\{\beta\})=\sum_{A\subset G}   q^{k(A)} \prod_{e\in A}\beta_e,
\ee
where $k(A)$ is the number of connected components of the subgraph $A$.

\bigskip

\begin{definition}
\label{defBR}
The Bollob\'as-Riordan polynomial of a ribbon graph $G$ is defined as
\begin{eqnarray}
R_G(x,y,z)
=\sum_{H\subset G}(x-1)^{r(G)-r(H)}y^{n(H)}z^{k(H)-F(H)+n(H)}.
\end{eqnarray}
\end{definition}

Note that we have denoted by $F(H)$ the number of components of the boundary of the respective subgraph $H$.

As the Tutte polynomial, the Bollob\'as-Riordan polynomial also obeys the deletion/contraction relation:

\begin{theorem}
\label{delcontrbollo}
Let G be a ribbon graph. One then has
\begin{equation}
\label{cdBR}
R_G=R_{G/e}+R_{G-e}
\end{equation}
for any regular edge $e$ of $G$. 
\end{theorem}

\begin{figure}
\begin{center}
\includegraphics[scale=0.6]{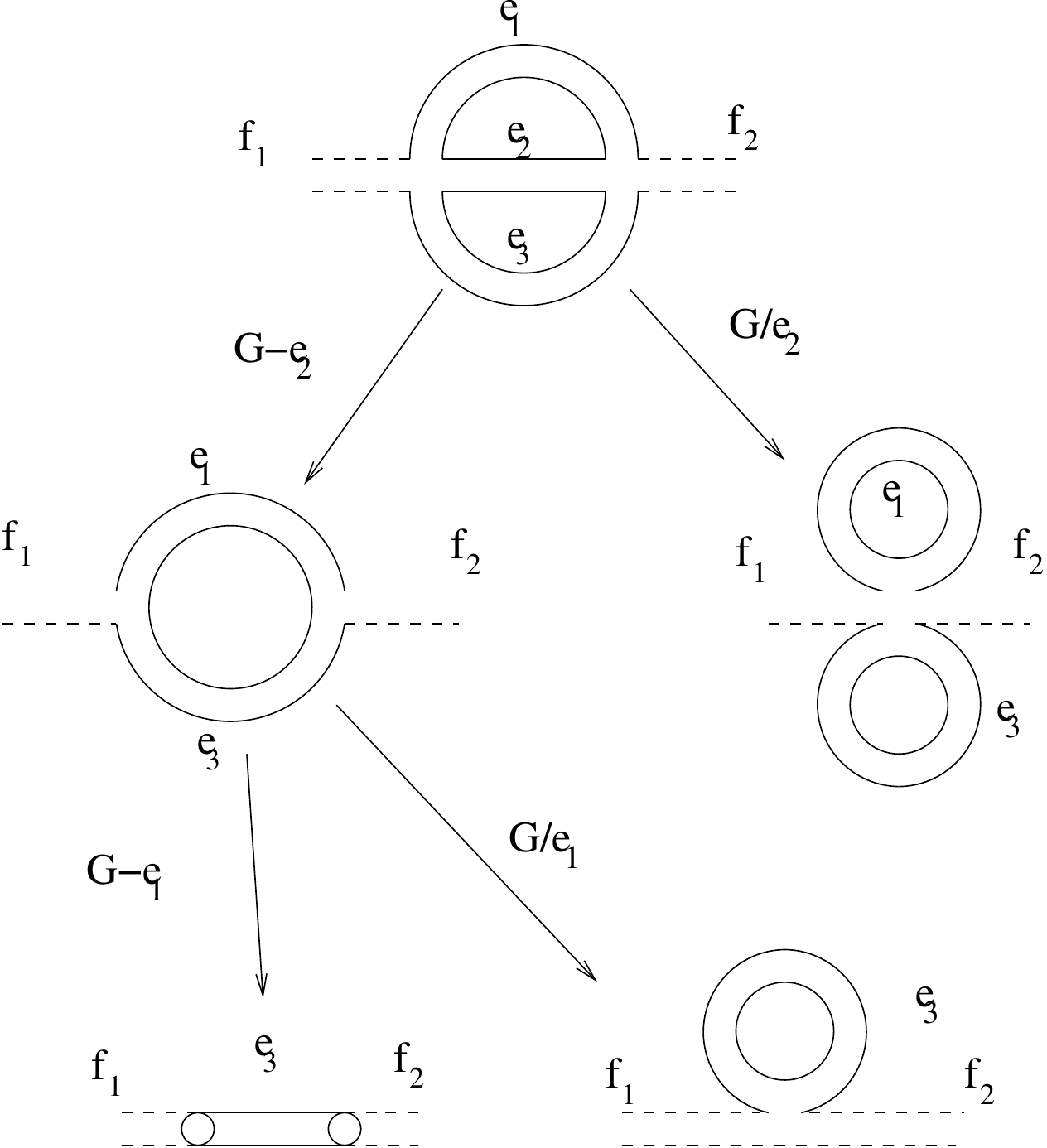}
\caption{The deletion/contraction of the regular edges of some ribbon graph $G$.}
\label{fig:r}
\end{center}
\end{figure}

To illustrate this deletion/contraction property for a ribbon graph, an explicit example is given in Fig. \ref{fig:r}. Note that flags (or external edges) have been this time added to the graph (the vertex valence being always equal to four - the $\Phi^4$ model).

\medskip

A multivariate version of the Bollob\'as-Riordan polynomial also exists in the literature:

\begin{eqnarray}
V_G(x,\{\beta\},z)
=\sum_{H\subset G}x^{k(H)}\left(\prod_{e\in H}\beta_e\right)z^{F(H)}.
\end{eqnarray}

\medskip

A signed version of the Bollob\'as-Riordan 
was also defined in
 \cite{signedBR}. Some partial duality was then defined in \cite{pd}; moreover, the properties of the multivariate version of this signed  Bollob\'as-Riordan polynomial were analyzed in \cite{fab} (namely, its invariance under the partial duality of \cite{pd} was proved). A further contribution was made in \cite{sur}, where a four-variable generalization of the Bollob\'as-Riordan polynomial for ribbon graphs was defined.

\section{Diagrammatics of GFT - bubbles and colorable graphs}
\label{GFT}
\renewcommand{\theequation}{\thesection.\arabic{equation}}
\setcounter{equation}{0}

In this section we present some notions of diagrammatics of GFT, emphasizing on the definitions of bubbles and of colorable graphs.

GFT is defined as a quantum field theory on group manifolds. Its associated Feynman graphs are, for the three-dimensional case, rank three tensor graphs. 

\begin{figure}
\begin{center}
\includegraphics[scale=0.3,angle=0]{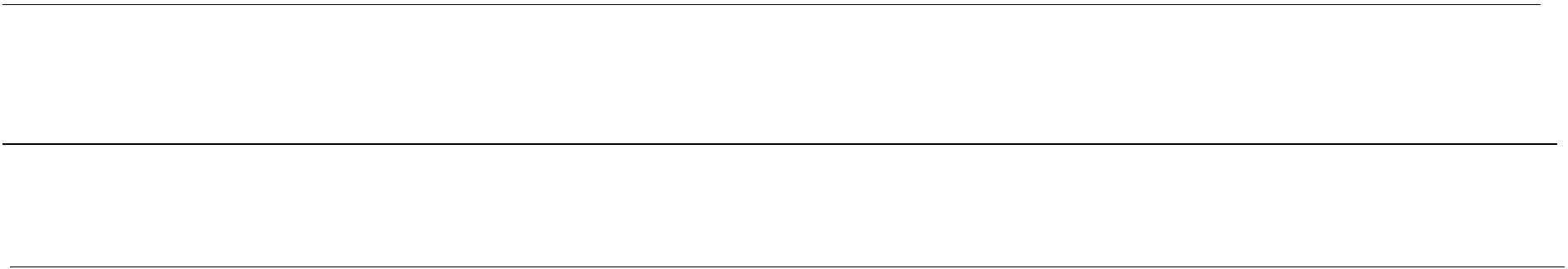}
\caption{The edge of the GFT tensor graphs.}
\label{fig:edge}
\end{center}
\end{figure}

In order to define such a scalar Feynman graph one needs a propagator (associated to the edge of the graph) and an interaction (associated to the graph vertex). The edge is simply given by Fig. \ref{fig:edge}.  Each such GFT tensor graph edge has three strands.

\medskip

The vertex of the three-dimensional GFT is a valence four vertex. 
Let us give some more explanation on this.

\begin{figure}
\begin{center}
\includegraphics[scale=0.4,angle=0]{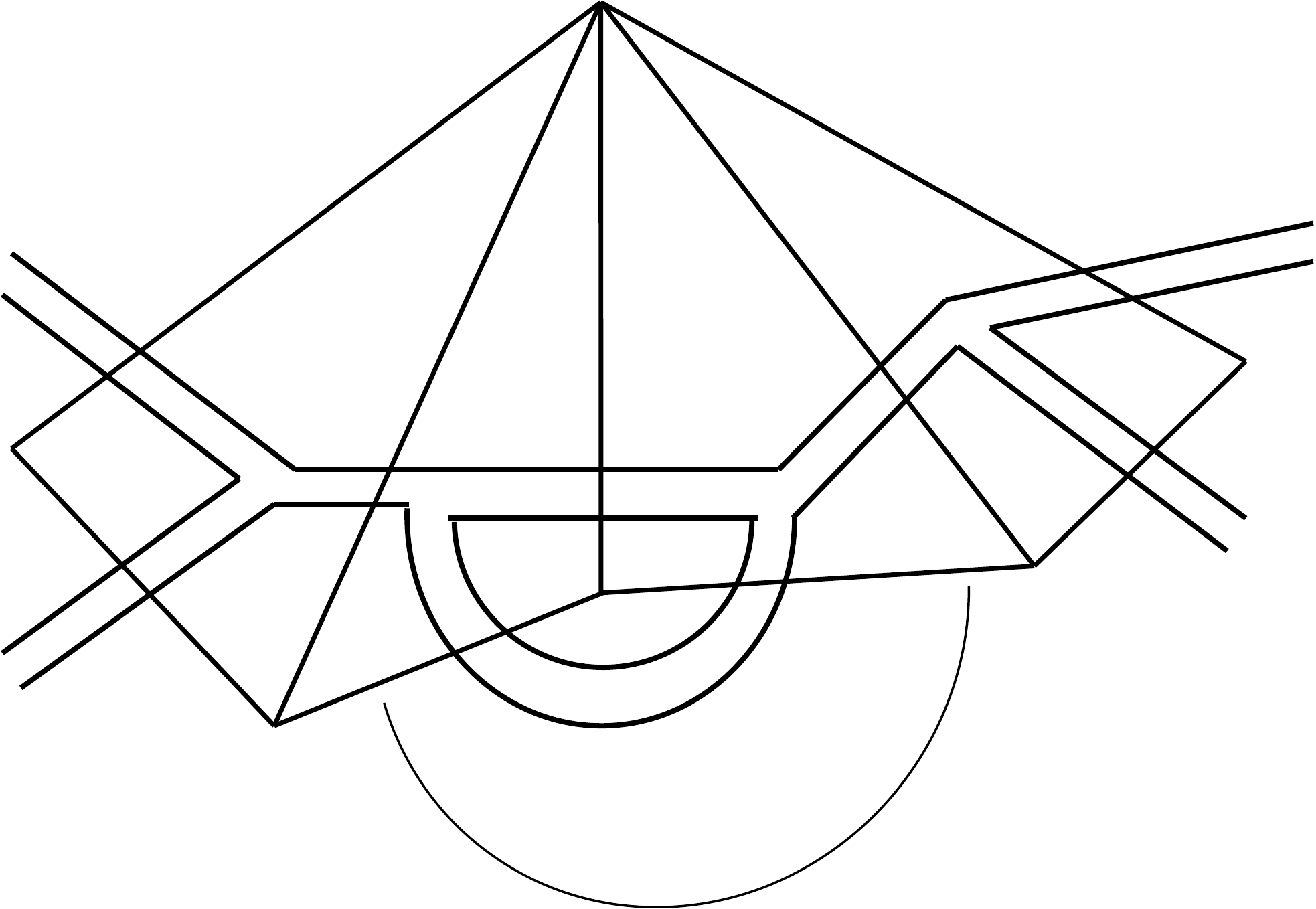}
\caption{Some triangularization of a two-dimensional surface and the dual point of view leading to valence three ribbon graphs.}
\label{fig:2d}
\end{center}
\end{figure}

\begin{figure}
\begin{center}
\includegraphics[scale=0.4,angle=0]{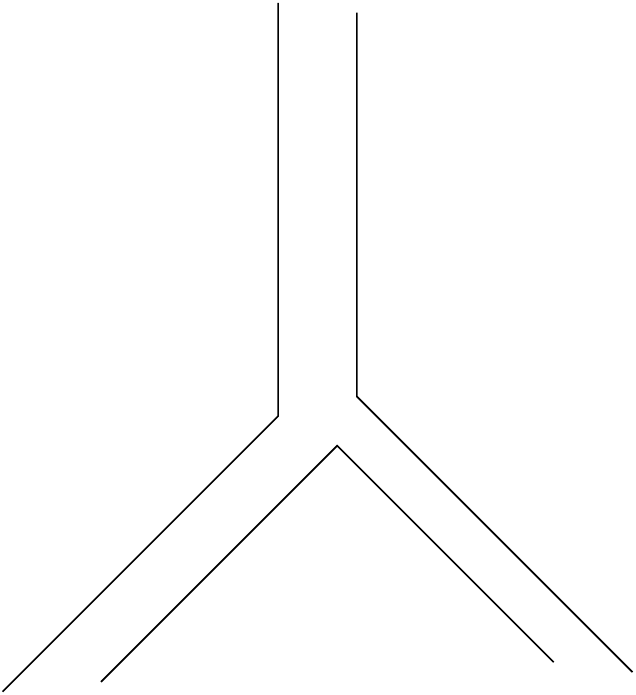}
\caption{The valence three vertex of two-dimensional GFT.}
\label{fig:vertex2d}
\end{center}
\end{figure}

As already mentioned in the introduction, GFT graphs are dual to triangularizations of space-time. 
Thus, in the simplest two-dimensional case, the building block of such a triangularization is of course a triangle. One can now take a dual point of view and encode the information coming from this triangularization in a ribbon graph picture (see for example the two-dimensional triangularization of Fig. \ref{fig:2d}). 
The appropriate valence of the associate ribbon graph vertex is thus three: the vertex (see Fig. \ref{fig:vertex2d}) corresponds to the triangle (the building block of the triangularization) and the three incoming/outgoing edges correspond to the three edges of the original triangle. 
The two strands of each edge correspond to the two vertices joined by the respective edge in the initial triangle. For the sake of completeness, let us also mention that this point of view relates to the one of two-dimensional quantum gravity matrix models \cite{david}.

\begin{figure}
\begin{center}
\includegraphics[scale=0.4,angle=0]{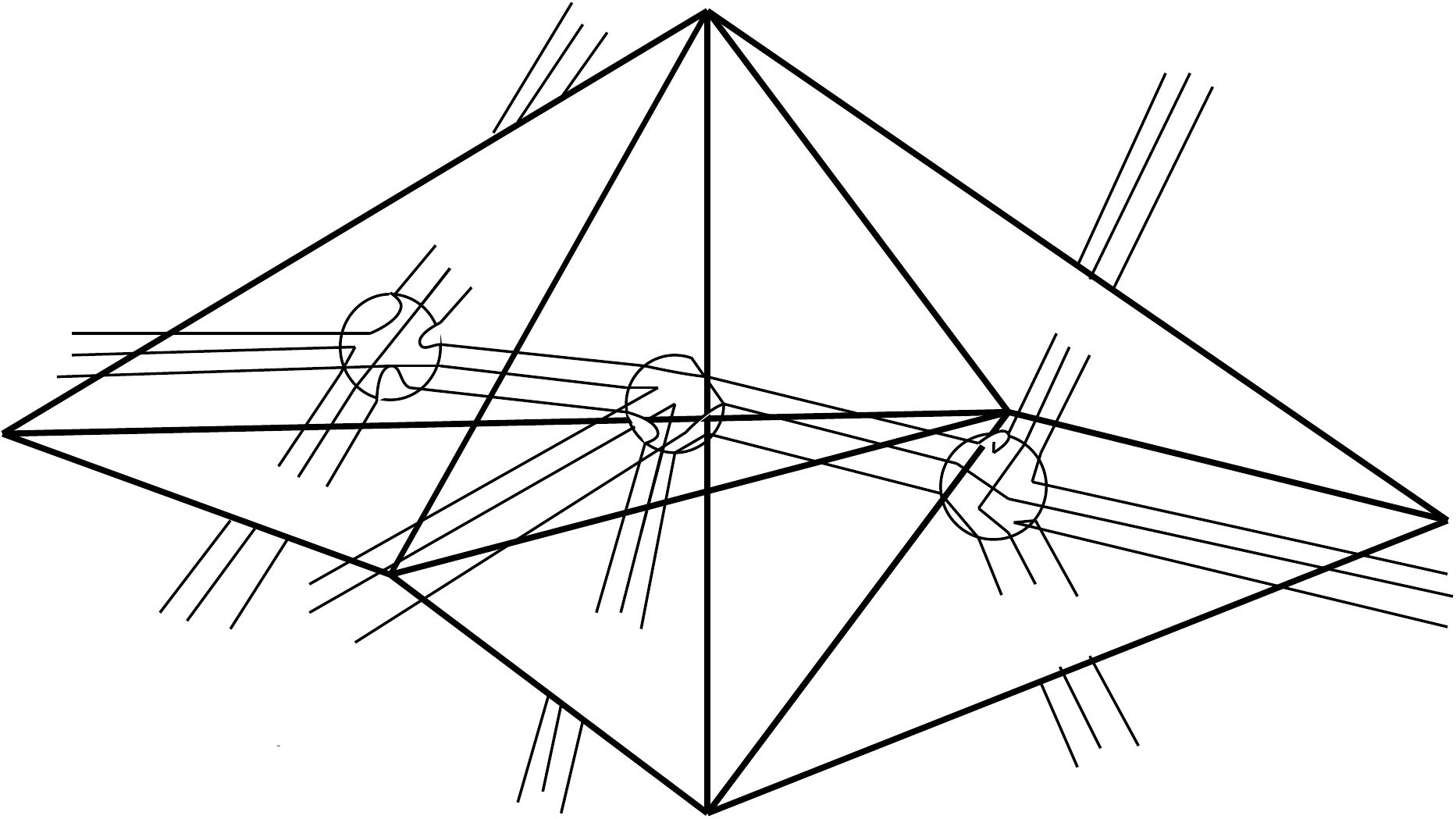}
\caption{Some triangularization of a three-dimensional surface and the dual point of view leading to valence four tensor graphs.}
\label{fig:t}
\end{center}
\end{figure}

The situation is similar for the three-dimensional case. The valence four vertex corresponds to a tetrahedron, the simplest building-block of a three-dimensional space-time. The four triangles constructing the tetrahedron correspond to the four edges intersecting at the respective vertex (see for example Fig. \ref{fig:t}). 
The three strands of such an edge correspond to the the three edges of the respective triangle (face of the tetrahedron), thus generalizing the image of the two-dimensional case. 
Let us recall at this point that more details on GFTs and their physical meaning are given for example in \cite{gft}.

\medskip

For the sake of simplicity, we analyze in this paper orientable three-dimensional tensor graphs. This means that the simplicial complexes dual to these graphs are orientable (the interested reader can turn for example to \cite{fgo} for more details on this point). Let us also mention that the concepts used in this paper do not seem to apply directly to general tensor graphs.

\begin{figure}
\begin{center}
\includegraphics[scale=0.4,angle=0]{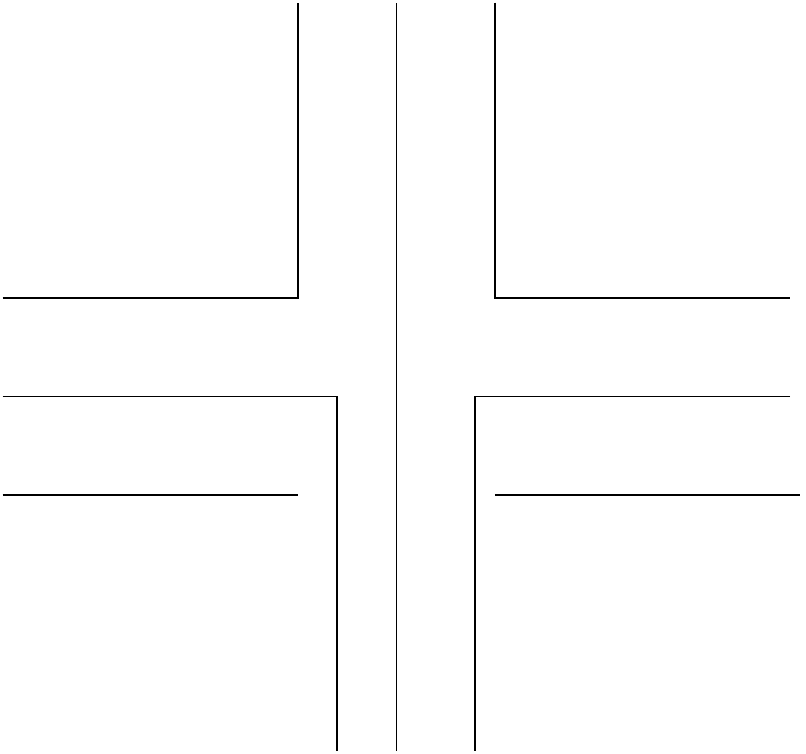}
\caption{The vertex of the GFT tensor graphs.}
\label{vertex}
\end{center}
\end{figure}

One can use several conventions for drawing the vertex of such a theory. For later convenience, we use the one of \cite{fgo}, given in Fig. \ref{vertex}. Note that other conventions do not lead to different results of this paper.

\begin{figure}
\begin{center}
\includegraphics[scale=0.35,angle=0]{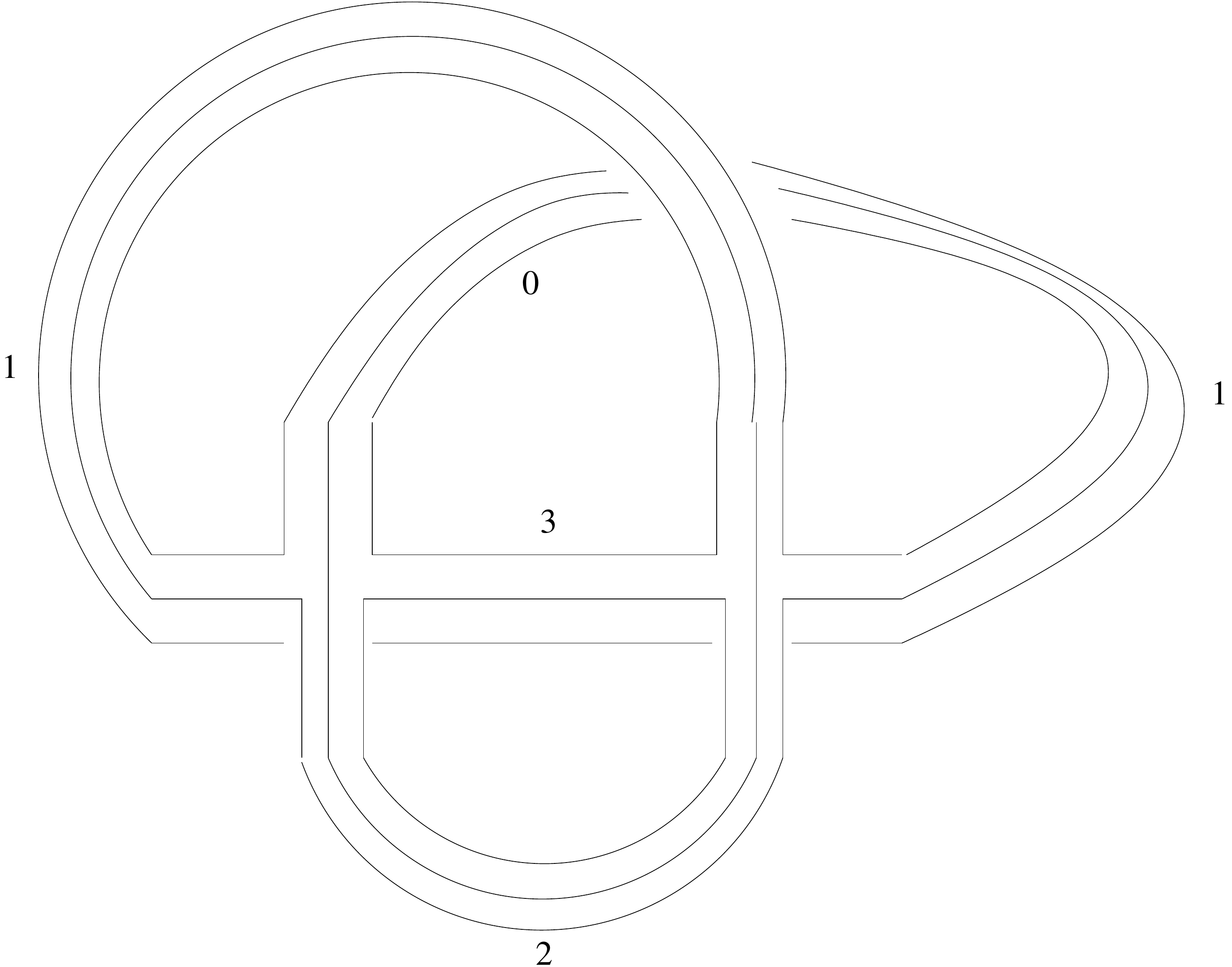}
\caption{An example of a GFT graph. We denote by $e_0,\ldots, e_3$ the four edges.}
\label{nc}
\end{center}
\end{figure}

Assembling together these two constituents, one is now able to build up tensor graphs. An example of such a graph is the one given in Fig. \ref{nc}.

\medskip

Let us emphasize that, as mentioned in the Introduction, these tensor graphs are a natural generalization of the ribbon graphs described by the Bollob\'as-Riordan polynomial; the strand structure described above makes these theories extremely rich from a combinatorial and topological point of view.

Thus, in addition to {\it faces} (closed circuits) which are 
natural to define 
for ribbon graphs, one can further generalize this topological notion for the case of tensor graphs:

\begin{definition}
A bubble is a closed three-dimensional region of the graph.
\end{definition}

A bubble can be represented as a ribbon graph, and thus, from a topological point of view, it is natural to define the its Euler characteristic $\chi$ (and hence its  genus).  One has
\beqa
\label{euler}
\chi=2-2g=V-E+F,
\eeqa
where $g$ is the bubble genus
 and, as before, $V$ is the total number of vertices, $E$ is the total number of internal edges and $F$ is the total number of faces. 
This relation is just a translation of the fact that such a ribbon graph can be drawn on a two-dimensional manifold (each graph defining the surface on which it is drawn), $g$ thus being the genus of the respective manifold.
One thus speaks of {\it planar} or {\it non-planar bubbles}.

\begin{figure}
\begin{center}
\includegraphics[scale=0.35,angle=0]{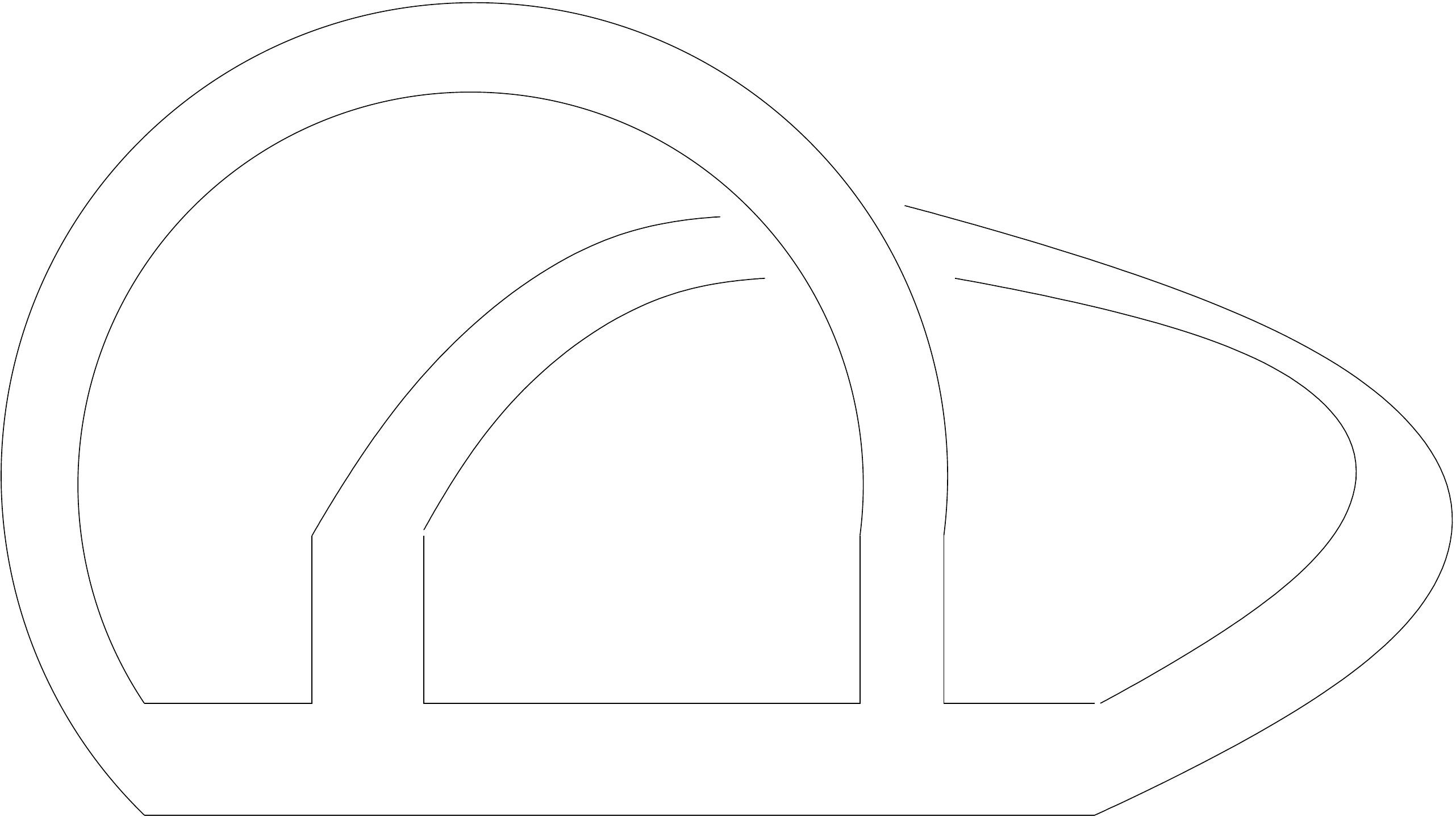}
\caption{One of the bubbles of the graph of Fig. \ref{nc}. It has two vertices, three edges and $1$ face and thus genus $1$; it is  non-planar.}
\label{bula1NC}
\end{center}
\end{figure}
\begin{figure}
\begin{center}
\includegraphics[scale=0.35,angle=0]{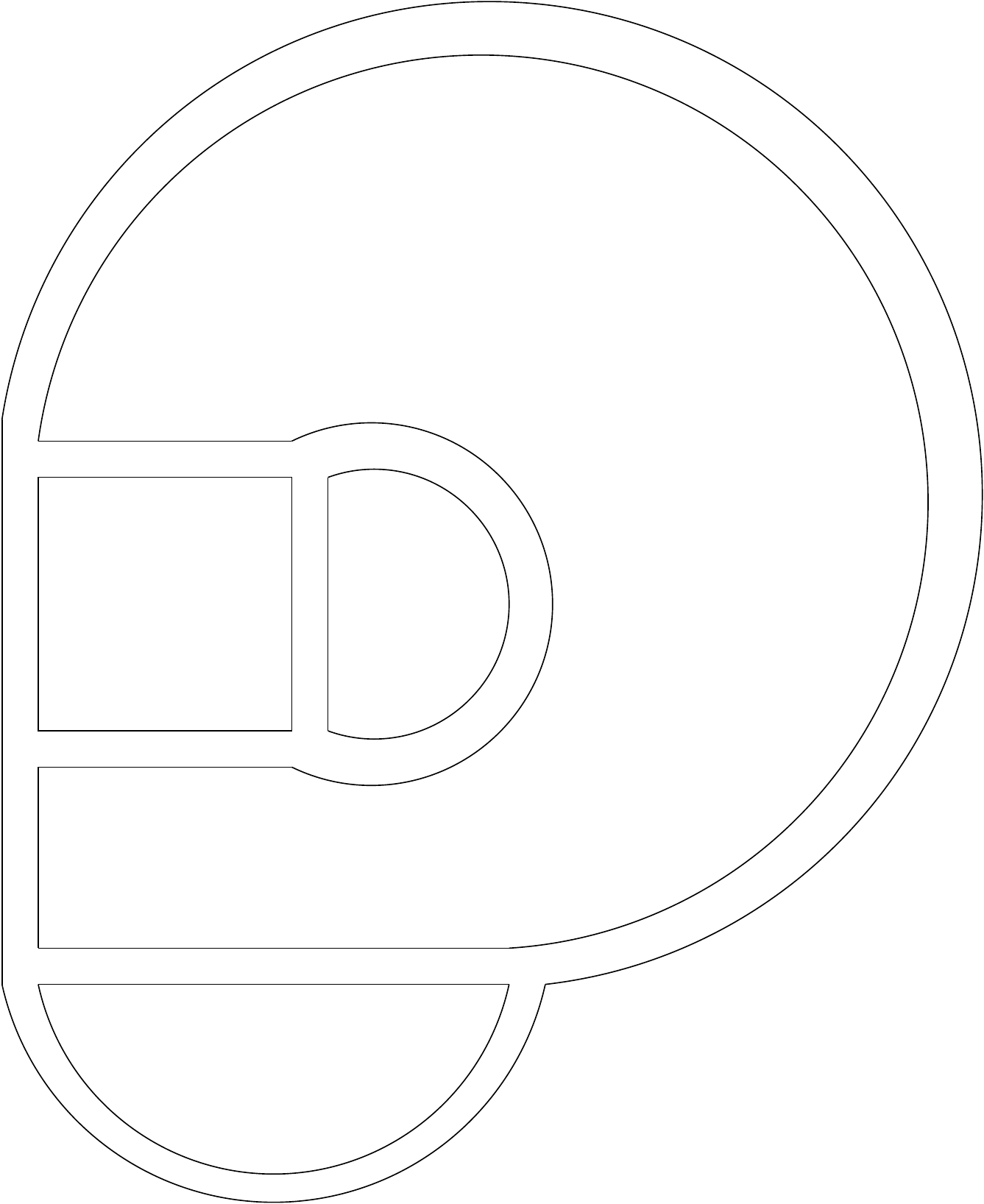}
\caption{The second of the bubbles of the graph of Fig. \ref{nc}. It has  genus $0$, being hence planar.}
\label{bula2NC}
\end{center}
\end{figure}

For the example of Fig. \ref{nc}, the bubbles are given in Fig. \ref{bula1NC} and \ref{bula2NC}. The first of them is non-planar (having a single face and thus, according to \eqref{euler}, genus $1$), while the second of them is planar.

Note that in \cite{fgo} a precise algorithm for characterizing these bubbles was given for an arbitrary GFT tensor graph. A cohomological definition of this notion of bubble was given in \cite{matteo}. 

As faces do for ribbon graphs, bubbles are known  to play a crucial r\^ole for these tensor graphs. Thus, if in QFT the divergences (which require renormalization) are related to cycles (known under the name of loops in mathematical physics), in GFT these divergences are related to this new notion of bubbles.

Furthermore, another important result which we will exploit in the sequel, is the following. As already mentioned above, three-dimensional tensor graphs correspond to triangularizations of  three-dimensional spaces. Furthermore, if all the bubbles of the graph are planar, the respective graph represents a {\it manifold}. If any of the graph's bubbles is non-planar, then the graph represents (at least) a {\it pseudo-manifold} (in the three-dimensional case we analyze here they only present isolated singularities) 
\cite{lost}.

We thus notice the particular importance of the sum of the genera of these bubbles. As already announced in the introduction, it is this variable which we will use to generalize the Bollob\'as-Riordan polynomial.

\bigskip

Before going further let us also present the notion of {\it colorable} tensor graphs \cite{color}. These graphs represent a class of graphs which have a simpler mathematical behavior. 

\begin{definition}
A colorable graph is a graph where each edge can be assigned a certain color belonging to the set $\{c_0,\ldots,c_3\}$, such that 
each of the four incoming/outgoing edges of an arbitrary vertex has a distinct color and 
at each vertex one has a clockwise or an anti-clockwise cyclic ordering for the colors $\{c_0,\ldots,c_3\}$ of the four incoming/outgoing edges.
\end{definition}

\begin{figure}
\begin{center}
\includegraphics[scale=0.35,angle=0]{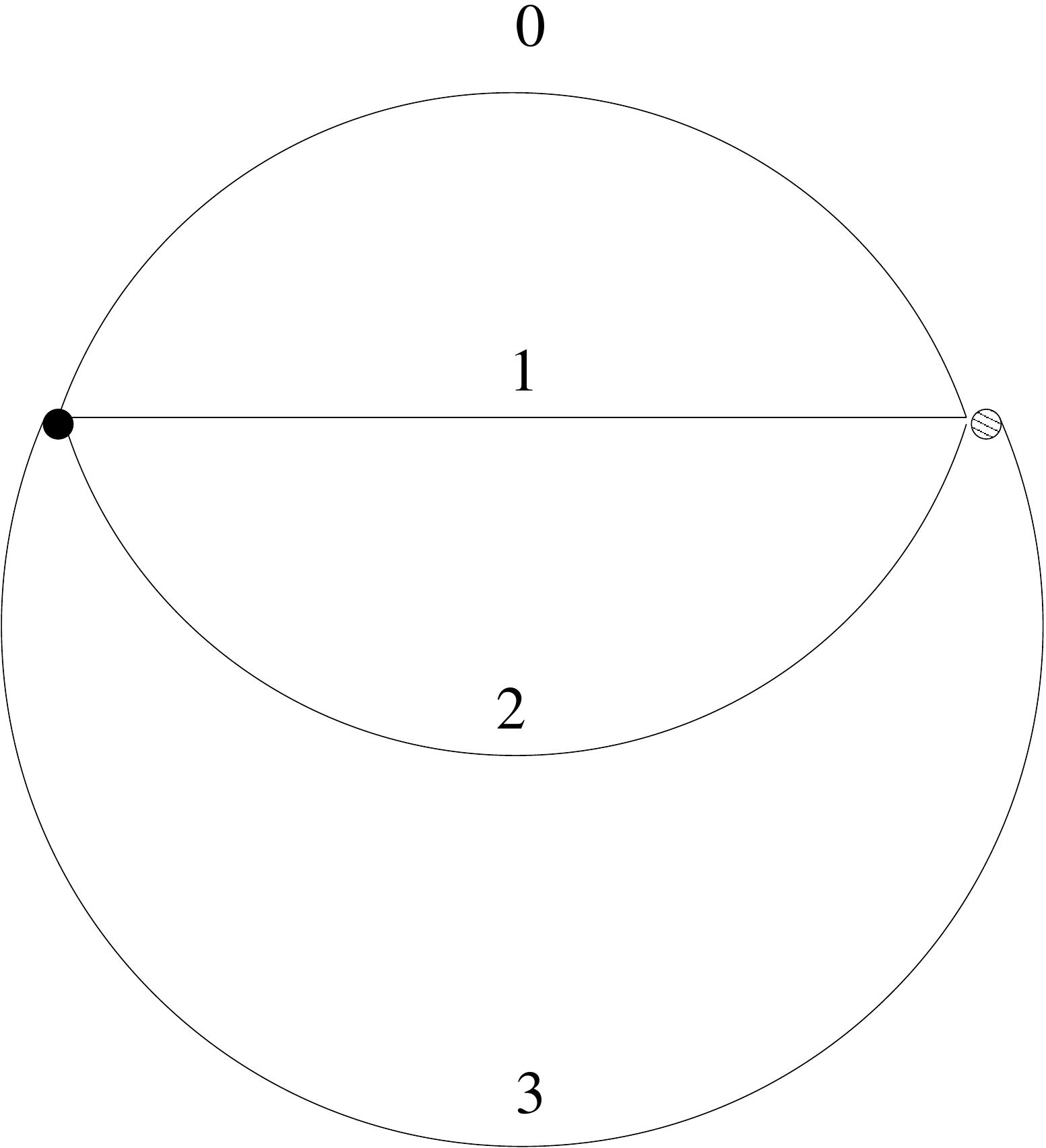}
\caption{An example of a colorable graph. The four colors are assigned with the numbers $0,\ldots,3$ associated to the edges. The graph has two vertices, a white and a black one. At the black vertex, the four colors $0,\ldots,3$ come in a clockwise ordering, while at the white vertex, one has an anti-clockwise ordering of the colors.}
\label{col}
\end{center}
\end{figure}

We are dealing here with only four possible colors, since we are interested in the three-dimensional case.
As a direct consequence of the definition above, one can conclude that a colorable graph has two types of vertices: a {\it white (or positive)} one (where the order is anti-clockwise) one and a {\it black (or negative) one} (where the order is clockwise). 
Furthermore, any edge of such a graph connects a positive to a negative vertex.
An example of such a colorable graph is the one of Fig. \ref{col}. Note that the graph of Fig. \ref{nc} is not colorable.

Furthermore, let us emphasize on the following. When drawing these colorable graph, one can drop the cumbersome stranded standard picture, because all the information is coded in the coloring. This can actually be done for any GFT tensor graph (colorable or not) which does not  allow twists in the edges (or vertices), see also the definition of jacket graphs below.

In QFT language the coloring presented here is related to the fact that one deals with a complex field $\Phi$. Thus, the two types of vertices correspond to two types of interaction: a $\Phi^4$ one and a $\bar \Phi^4$ one, where we have denoted by $\bar \Phi$ the complex conjugate of the field $\Phi$.


\medskip

\begin{figure}
\begin{center}
\includegraphics[scale=0.35,angle=0]{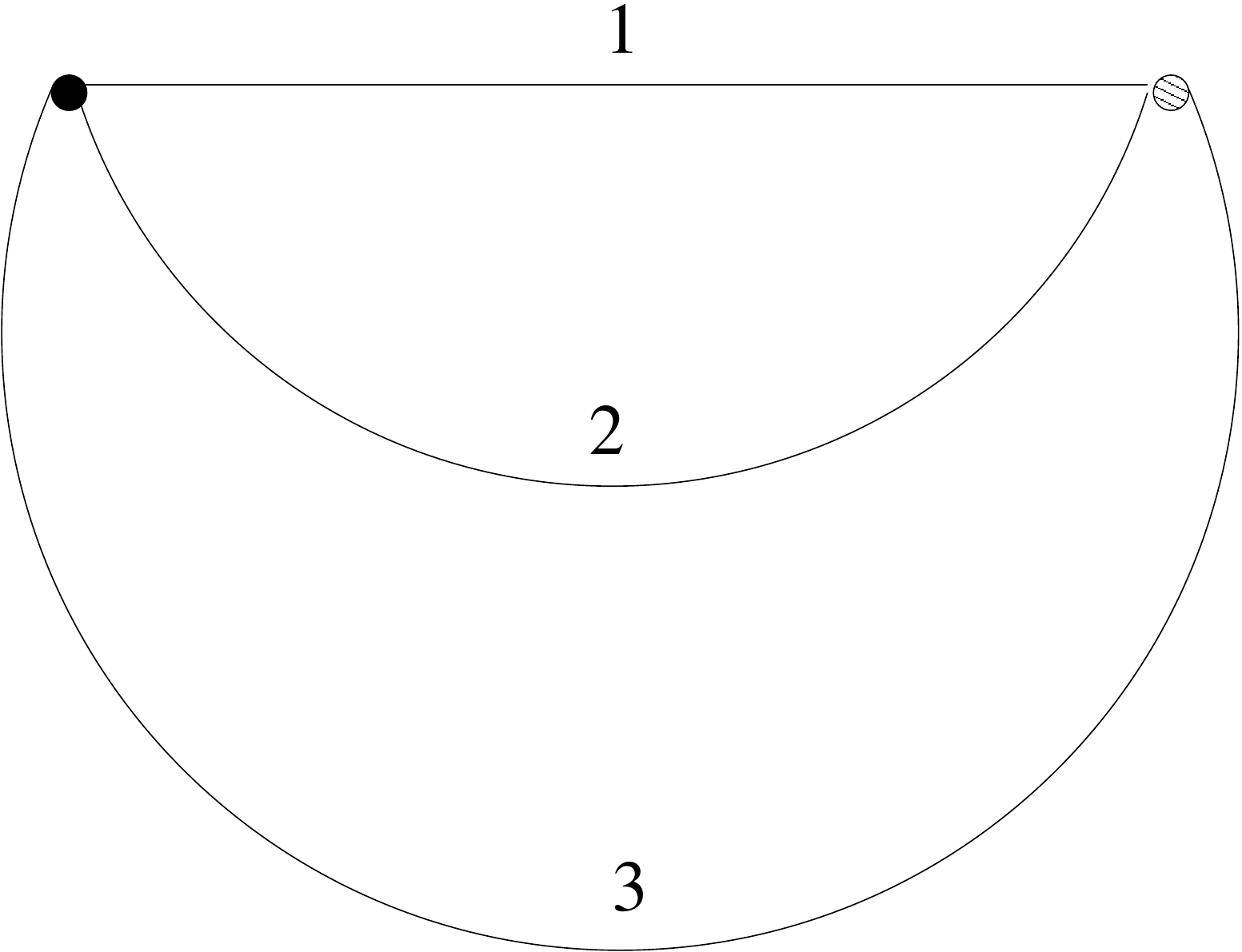}
\caption{The first of the bubbles of the graph \ref{col}; it is obtained by deleting the edges of color $0$ of the initial graph.}
\label{b1}
\end{center}
\end{figure}

\begin{figure}
\begin{center}
\includegraphics[scale=0.35,angle=0]{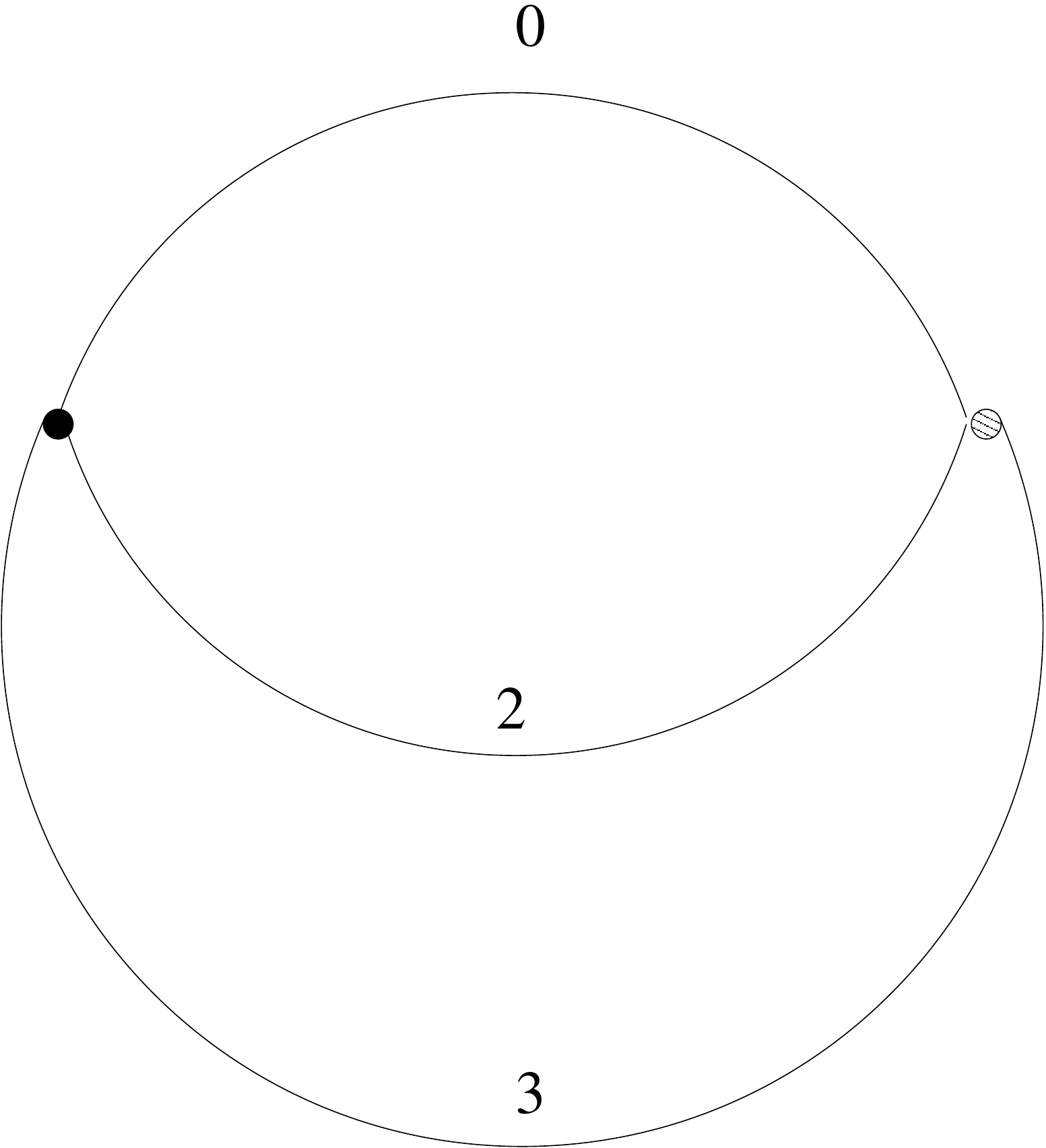}
\caption{The second of the bubbles of the graph \ref{col}; it is obtained by deleting the edges of color $1$ of the initial graph.}
\label{b2}
\end{center}
\end{figure}

\begin{figure}
\begin{center}
\includegraphics[scale=0.35,angle=0]{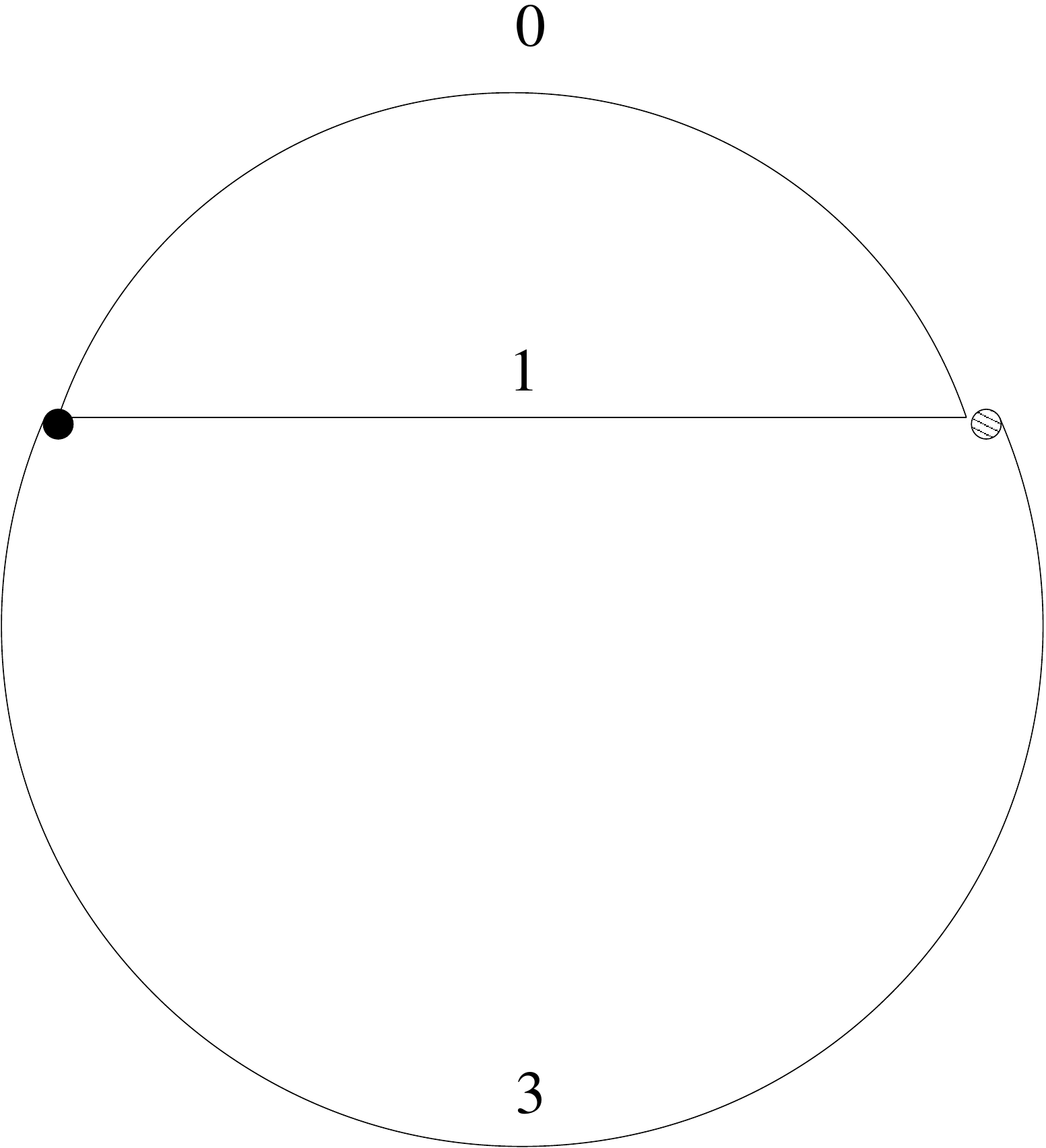}
\caption{The third of the bubbles of the graph \ref{col}; it is obtained by deleting the edges of color $2$ of the initial graph.}
\label{b3}
\end{center}
\end{figure}

\begin{figure}
\begin{center}
\includegraphics[scale=0.35,angle=0]{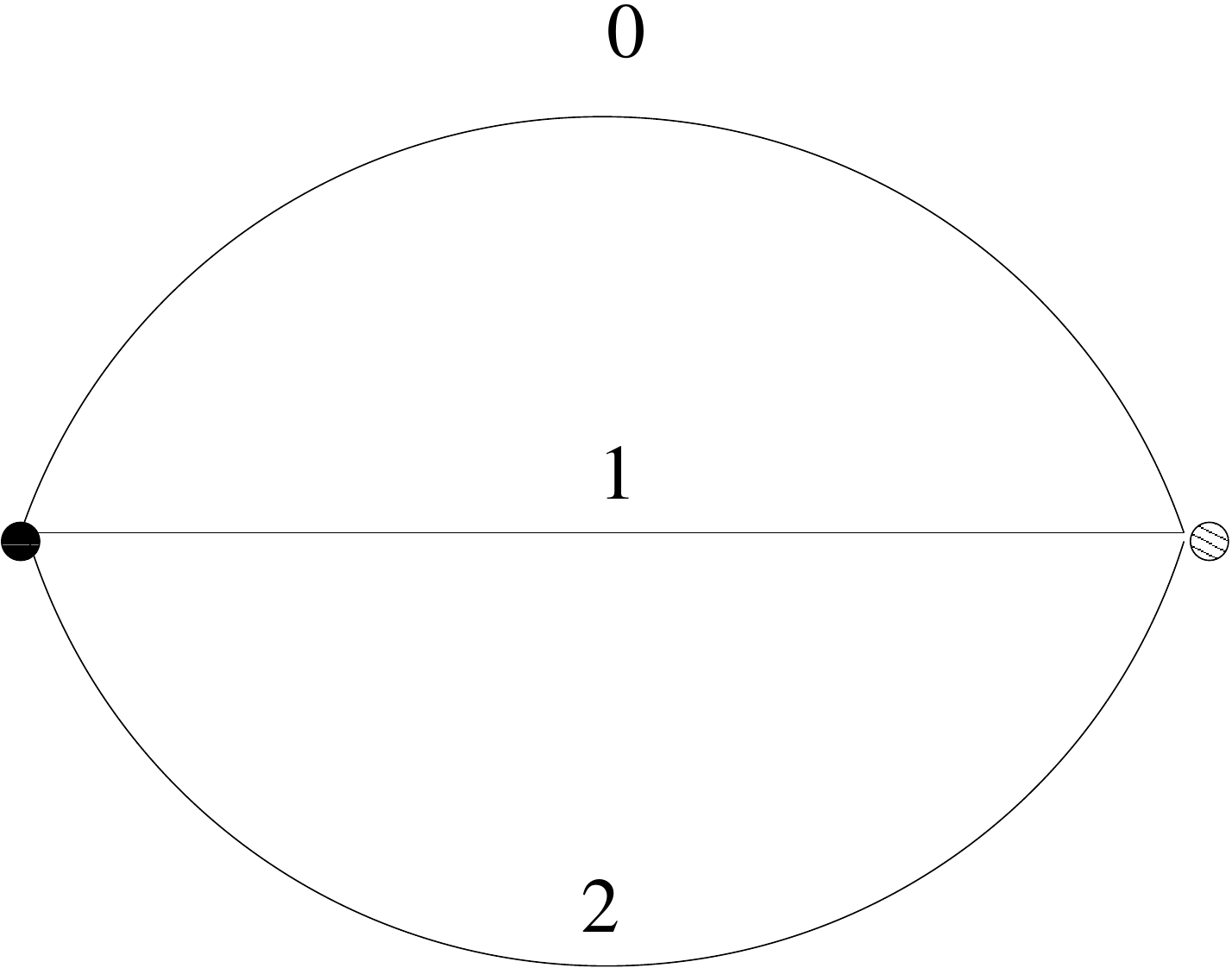}
\caption{The last of the bubbles of the graph \ref{col}; it is obtained by deleting the edges of color $3$ of the initial graph.}
\label{b4}
\end{center}
\end{figure}

As expected, the topological notion of bubbles is defined in an easier way for these colorable graphs. Thus, a bubble can now be defined as a connected component of a subgraph formed only of edges of colors belonging to a subset of cardinal three of the initial set of colors (see again \cite{lost, color}).
For the colorable graph of Fig. \ref{col}, the bubbles are drawn in Fig. \ref{b1}, \ref{b2}, \ref{b3} and \ref{b4}.

\medskip

Before ending this section, let us remark that the GFT tensor graphs presented in this paper (colorable or not) can be graphically represented as ribbon graphs, if one erases one of the 
strands of each edge. With the convention of Fig. \ref{vertex}, this can be done by erasing the middle strand for the vertical edge and respectively the lower strand for the horizontal edge. 
The set of graphs obtained in such way is in bijection with the initial set of GFT tensor graphs, because of the fact that one does not allow permutations (or twists) in the edge of the model (or in the vertex).  The graph thus obtained is referred to as a {\it jacket graph} in \cite{gft-regain-2}. 


One can thus raise the legitimate question of the r\^ole played by this middle strand. One of the answers is that without the middle strand one cannot define the bubbles of the graphs and obviously one does not have the full topological richness of this class of graphs. Moreover, let us also remark that the non-commutative scalar quantum field theories which are known to be renormalizable on Moyal space, also involve ribbon graphs which do not allow twists at the level of the graph edges (or vertices) (see for example \cite{io-NCQFT} and references within).

Furthermore, the strand structure also plays a crucial r\^ole when conserving, in GFT analytic calculations,  physical quantities like the ``momentum''.

\section{Tensor graph contraction; the Gur\u au polynomial}
\label{Gurau}
\renewcommand{\theequation}{\thesection.\arabic{equation}}
\setcounter{equation}{0}

In this section we recall the operation of contraction for tensor graphs introduced in \cite{gp}; for the sake of completeness, we also give the definition of the Gur\u au polynomial \cite{gp}.

Note that, when contracting an usual graph, the valence of its vertices can increase. This is not a problem for a ribbon graph (faces can still be defined in a straightforward manner even for such a contracted ribbon graph) but for the tensor graphs presented in this paper, one has to define a way of matching the middle strands of the respective edges. Furthermore it is not clear if the notion of bubble can still be defined for such contracted graphs.

In the case of colorable graphs, the situation is even more dramatic, because the contracted graph is not a colorable graph anymore.

To circumvent this problem, a different way of contracting tensor graphs was proposed in \cite{gp}. The definition was given for a colorable graph, but we propose here to extend to it to an arbitrary GFT graph.

Thus, when an edge is contracted, it just becomes {\it passive}. The set of edges of the graph $G$ is now formed of the disjoint reunion of two subsets, the active ones, denoted by ${\cal L}_1(G)$, and the passive ones, denoted by ${\cal L}_2(G)$. Hence the tensor graph contraction transforms an active edge into a passive one (instead of  identifying its two vertices, as it is happening in the Tutte or the Bollob\'as-Riordan cases). The deletion (which is defined as usually) and the contraction operations are then allowed only for active lines.

The Gur\u au polynomial associated to a colorable graph $G$ is then defined as
\beqa
\label{gurau}
{\mathcal G}_G \left(\{\beta\},\{x\},\{y\}\right)=
\sum_{H\in G,\ {\cal L}_2(H)}\left( \prod_{e\in{\cal L}_1(H)}\beta_e\right)
\prod_{p=0}^{n+1}x_p^{|{\mathcal B}^p|}\prod_{p=0}^{n}y_p^{|{\mathcal B}^p_\partial|}.
\eeqa
In this formula, $n+1$ is the number of colors of the graph. 
As already explained in the previous section, the value of interest for us is $n=3$ (since we deal in this paper with the three-dimensional case).
The set of variables $\beta$ is associated to the edges of the graph, while the set of variables $x$ is associated to the topology of the graph: one keeps track of the number of vertices, edges, faces, bubbles and finally connected components of the graph. The set of variables $y$ plays a similar r\^ole, but is associated to some boundary graph defined from the initial tensor graph (see again \cite{gp} for details).

Note that this polynomial is a multivariate one, in the sense that one has a distinct variable for each edge. Furthermore, it is  proved that it also satisfies the deletion/contraction relation (see Lemma $4$ of \cite{gp}).

\medskip

Let us end this section by stating that recently, a generalization of the Bollob\'as-Riordan polynomial was proposed in \cite{surpriza} for triangulations. The point of view adopted is completely different of the one here (or the one in \cite{gp}), no GFT tensor graphs being used.

\section{Definition of the ${\mathcal T}$ polynomial; deletion/contraction and multivariate version}
\label{T}
\renewcommand{\theequation}{\thesection.\arabic{equation}}
\setcounter{equation}{0}

In this section we give the definition of a polynomial different of the Gur\u au polynomial. We then show that this new polynomial also satisfies the deletion/contraction relation. Finally, we define a multivariate and a hypervariate version of this polynomial.

The main idea behind this proposition is, as already announced before, to exploit the importance of the sum of genera of the bubbles of the graph. As we have already mentioned in section \ref{GFT}, if this sum is vanishing, then the respective graph is dual to a manifold, otherwise it is dual to a pseudo-manifold or to an even more singular surface.

Let us recall that the Bollob\'as-Riordan polynomial generalized the Tutte polynomial by adding a supplementary variable $z$ in order to keep track of the increase of topological data (see section \ref{TBR}). We thus propose to keep the same scheme by adding one more variable, $t$ to the set of variables of the Bollob\'as-Riordan polynomial. This variable will keep track of this sum of genera of the bubbles, telling thus in a natural manner weather or not the respective graph is dual to a manifold or not. 

As already mention in the section \ref{GFT}, to the GFT tensor graphs studied here one can associate in an unique way ribbon jacket graphs. This allows to define the rank, nullity, genus, number of faces {\it etc.} of the respective graph in a straightforward manner.

One has:

\begin{definition}
\label{def}
If $G$ is a GFT tensor graph, then 
\beqa
\label{pol}
{\mathcal T}_G(x,y,z,t)=\sum_{H\subset G} (x-1)^{r(G)-r(H)}y^{n(H)}z^{k(H)-F(H)+n(H)}t^{2\sum_{{\cal B}(H)} g_{{\cal B} (H)}},
\eeqa
where we have denoted by ${\cal B}(H)$ the bubbles of the GFT tensor subgraph $H$. 
\end{definition}

As in the definition of the Bollob\'as-Riordan polynomial (see Definition \ref{defBR}), the notation $F(H)$ above stands for the number of components of the boundary of the respective subgraph $H$ (seen as a ribbon jacket graph).

Note that we have defined this graph  polynomial for an arbitrary GFT tensor graph, colorable or not, since one can define the bubbles (and hence calculate their genera) not only for colorable graphs (see again \cite{fgo}). Nevertheless, one can obviously define the polynomial \eqref{pol} only for colorable graphs, where the bubbles can be defined in terms of colors, without going through the strand structure of the respective graph (see section \ref{GFT}).

The factor two in the last exponent of \eqref{pol} is chosen for similitude with the Euler characteristics formula.

One has
\beqa
\label{comp}
{\mathcal T}_G(x,y,z,1)&=&R_G(x,y,z),\nonumber\\
{\mathcal T}_G(x,y,1,1)&=&T_G(x,y),
\eeqa
where, by an abuse of notation, we have denoted again with $G$ the (ribbon) graph associated to the tensor graph $G$ (from the left hand side of the equations).

\medskip

We now prove that this polynomial satisfies the deletion/contraction relation. Note that we use the contraction of tensor graphs described in the previous section, operation which extends naturally for an arbitrary GFT tensor graph, not only for a colorable one.

One has:

\begin{theorem}
\label{teorema}
If $G$ is a GFT tensor graph, then
\beqa
\label{th}
{\cal T}_G(x,y,z,t)={\mathcal T}_{G-e}(x,y,z,t)+{\mathcal T}_{G/e}(x,y,z,t),
\eeqa
if $e$ is an active regular edge.
\end{theorem}

{\it Proof.} When choosing the particular edge $e$ and writing the sum on the right hand side 
\eqref{th}, one has two types of terms. One involving the passive line $e$ (coming from the contribution of ${\mathcal T}_{G/e}(x,y,z,t)$) and one not involving it (coming from the contribution of ${\mathcal T}_{G-e}(x,y,z,t)$).

Let us now prove that all the contributions of subsets of edges of $G$ (from the left hand side) of \eqref{th} are present also on the right hand side (and no other contributions appear).

The contribution coming from all the edges of the initial graph $G$ comes on the right hand side from ${\mathcal T}_{G/e}(x,y,z,t)$ (when one considers all the edges, the active and the passive one). The contribution coming from the set of all the edges except one come in the right hand side in the following way. The contribution without the special edge $e$ comes from ${\mathcal T}_{G-e}(x,y,z,t)$. The contributions coming from the set of all the edges except one (different of $e$ this time) come from ${\mathcal T}_{G/e}(x,y,z,t)$.

The situation is analogous for all subsets of edges of the initial graph. The empty set contribution comes from ${\mathcal T}_{G-e}(x,y,z,t)$ (${\mathcal T}_{G/e}(x,y,z,t)$ having to keep track all the time of the passive line $e$ does not present such a contribution).

Note that no supplementary contributions appear on the right hand side with respect to the left hand side of \eqref{th}. Obviously, there is no mixing between the various contributions of different cardinal of subsets of edges. This completes the proof. (QED)

\bigskip

We now define a multivariate version of the polynomial \eqref{def}, generalizing in a natural manner the multivariate Bollob\'as-Riordan polynomial.

\begin{definition}
\label{def2}
If $G$ is a GFT tensor graph, then 
\beqa
\label{pol2}
{\mathcal T}_G(x,\{\beta\},z,t)=\sum_{H\subset E} x^{k(H)}\left(\prod_{e\in H}\beta_e\right)z^{F(H)}t^{2\sum_{{\cal B}(H)} g_{{\cal B} (H)}},
\eeqa
where we have denoted by ${\cal B}(H)$ the bubbles of the GFT tensor subgraph $H$.
\end{definition}

Note that the Gur\u au polynomial \eqref{gurau} is in this sense already a multivariate polynomial (being also naturally related to the multivariate Bollob\'as-Riordan and Tutte polynomials, as shown in \cite{gp}).


Furthermore, the polynomial \eqref{pol2} can be further generalized to a {\it hypervariate version}, if we make use of a set $\{\gamma_b\}$ of variables, one for each bubble of the graph, instead of a single variable $t$:

\begin{definition}
\label{def3}
If $G$ is a GFT tensor graph, then 
\beqa
\label{pol3}
{\mathcal T}_G(x,\{\beta\},z,\{\gamma\})=\sum_{H\subset G} x^{k(H)}\left(\prod_{e\in H}\beta_e\right)z^{F(H)}\prod_b \gamma_b^{2\sum_{{\cal B}(H)} g_{{\cal B} (H)}}.
\eeqa
\end{definition}

The multivariate and hypervariate versions defined above also satisfy the deletion/contraction relation.

\bigskip

Let us now comment on the universality of the polynomial introduced here. 
For the Bollob\'as-Riordan polynomial, this is proved in \cite{br} or, for the multivariate case in \cite{multibr}. A restatement of this universality property was proved in a different form in \cite{j}.
 
Such a proof is mainly made possible through a classification of the terminal forms of the ribbon deletion/contraction procedure (referred to as {\it chord diagrams} in \cite{br} or \cite{j} or as {\it rosettes} in the non-commutative QFT literature (see for example \cite{io-BR, io-NCQFT, KRVT, roz} and references within).

For the graph polynomial described in this paper, the situation is similar to the one described in \cite{gp} in the sense that, in order to have such a property, one needs to start with a classification of three-dimensional piecewise linear manifolds and pseudo-manifolds. Such a result is nowadays an open question in algebraic topology (see again \cite{gp} for a more detailed discussion).

\medskip

Let us also address the issue of the relation of the Gur\u au polynomial with the proposed polynomial \eqref{pol}, when restricting ourselves to colorable graphs. As already mentioned in the previous section, the Gur\u au polynomial depends on a set of variable (the $\beta$'s) corresponding to the edges of the graphs and on two more sets of variables, the first of them (the four variables $x$) associated to the number of vertices, edges, faces and bubbles, the second of these sets (the three variables $y$) being used in a similar way (but with respect to the topology of some boundary graph). This is somehow a completely different ``philosophy''  that the one used in this paper. Thus, the polynomial \eqref{pol} only needs four variable and even if one considers the multivariate version \eqref{pol2} it does not seem natural to look for a relation between these polynomials. For example, in the definition of the polynomial \eqref{pol} (or its multivariate version \eqref{pol2}), we do not take into account the topology of some boundary graph, as it was done in \cite{gp} through the use of the variables $y$. Furthermore, if the Gur\u au polynomial needs a set of variables for each edge and then four more variables $x$ to encode the topology of the three-dimensional colored GFT graph (one for number of vertices, one for the number of edges, one for the number of faces and one for the number of bubbles, as already explained above), we only need a total of four such variables in our approach, this being (another) built-in fundamental difference between the two approaches.

All this becomes more transparent when comparing the differences of the way the Bollob\'as-Riordan polynomial is obtained from the polynomial \eqref{pol} (see \eqref{comp}) with the way the multivariate Bollob\'as-Riordan polynomial is obtained from the Gur\u au polynomial:
\beqa
{\mathcal G}_G(\{\beta\}, \{1,1,z,x\},\{1,1,1\})=V_G(x,\{\beta\},z).
\eeqa

Let us address one more point with respect to this relation. The Euler characteristic of the simplicial complex dual to the orientable tensor graphs we deal with in this paper is equal to the sum of the genera of the bubbles. This characteristic appears also when one performs a particular rescaling of the variables of the polynomial of \cite{gp} (see again \cite{gp} for details). Nevertheless, it does not seem to us that there is deeper reason behind this arithmetic fact, since the Euler characteristic appears in the definition of the polynomial \eqref{pol}, which is obviously not the case for the polynomial  of \cite{gp}. No direct relation (like the one relating the Tutte to the Bollob\'as-Riordan polynomial or the Tutte (or the Bollob\'as-Riordan polynomial) to the polynomials of the parametric representation \cite{io-BR}) between these two polynomials can be obtained, for the reasons listed above.

\bigskip

Before ending this section, we stress that these definitions and properties also hold for graphs with { external edges}. As already mentioned in section \ref{TBR}, from a graph theoretical point of view, these external edges can be seen as {flags} attached to the vertices (taking care that the valence of the vertices remains always equal to four). Each such flag is attached to a single vertex. This notion is of particular importance in theoretical physics, because it is the external vertices which are directly connected to the various numbers measured in experiments.

\section{An example of a non-colorable graph}
\label{ex}
\renewcommand{\theequation}{\thesection.\arabic{equation}}
\setcounter{equation}{0}

Let us give in this section an example of a non-colorable 
GFT tensor graph illustrating the definition of the polynomial \eqref{pol} and also Theorem \ref{teorema}. 
We analyze here the graph of Fig. \ref{nc}. As already stated in section \ref{GFT}, the two bubbles of this graphs (see Fig. \ref{bula1NC} and respectively \ref{bula2NC}) have genera zero and respectively one. The definition of the polynomial \eqref{pol} thus allows to have a non-trivial dependence of the polynomial on the supplementary variable $t$. This comes from the contribution of the whole tensor graph $G$. 

One has
\beqa
\label{tot}
{\mathcal T}_G(x,y,z,t)=4+6y+(x-1)+2y^2z^2+2y^2+y^3z^2t^2.
\eeqa

\begin{figure}
\begin{center}
\includegraphics[scale=0.3,angle=0]{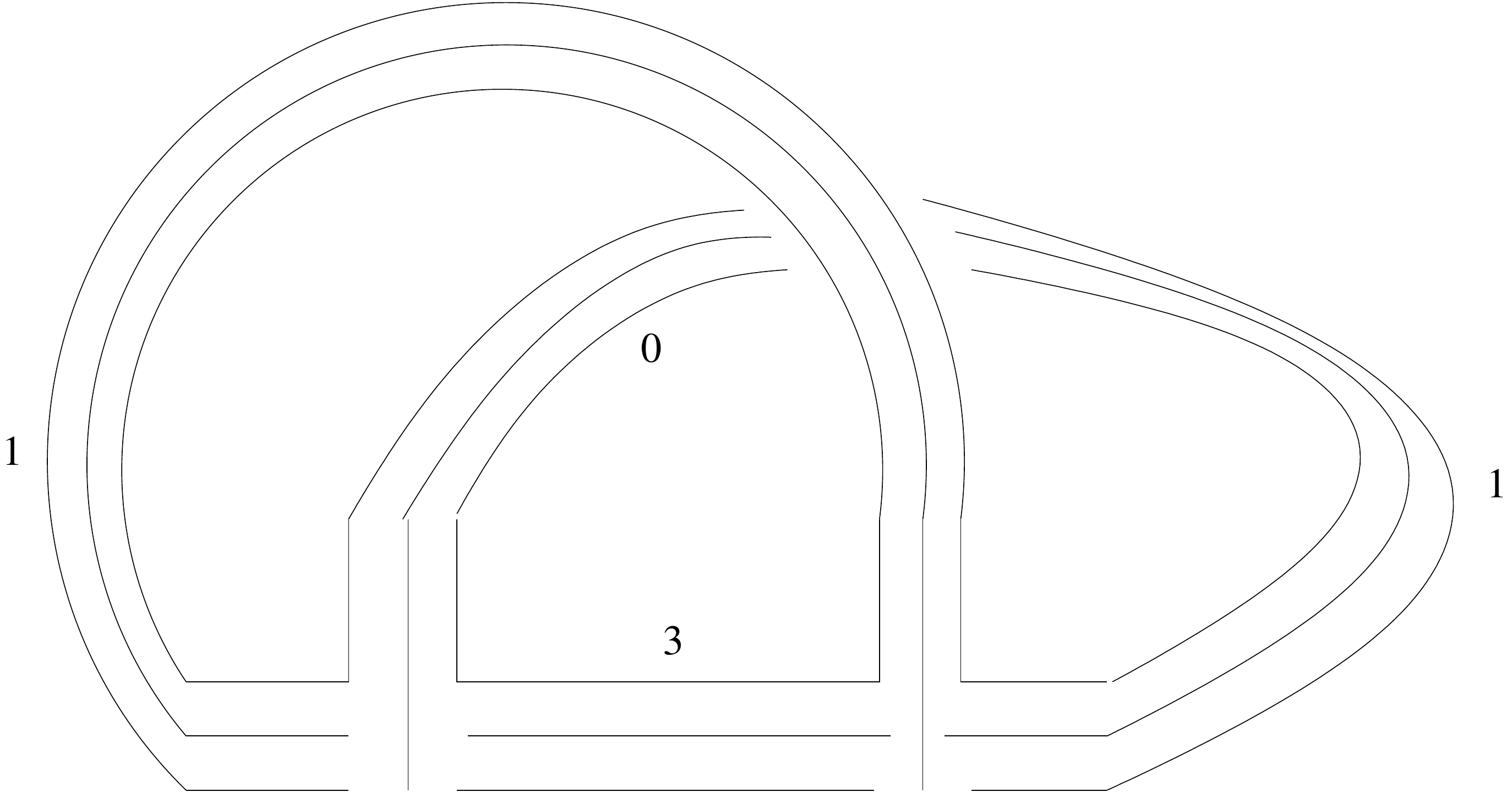}
\caption{The graph $G-e_2$ obtained from the graph $G$ of Fig. \ref{nc} by deleting the edge $e_2$.}
\label{G-e2}
\end{center}
\end{figure}

Let us choose the edge $e_2$ for the deletion and contractions operations. The graph $G-e_2$ is represented in Fig. \ref{G-e2} and leads to
\beqa
\label{con}
{\mathcal T}_{G-e_2}(x,y,z,t)=3+(x-1)+3y+y^2z^2.
\eeqa
Note that this graph has no bubbles; thus its associated polynomial will have trivial dependence in the supplementary variable $t$. 

\begin{figure}
\begin{center}
\includegraphics[scale=0.3,angle=0]{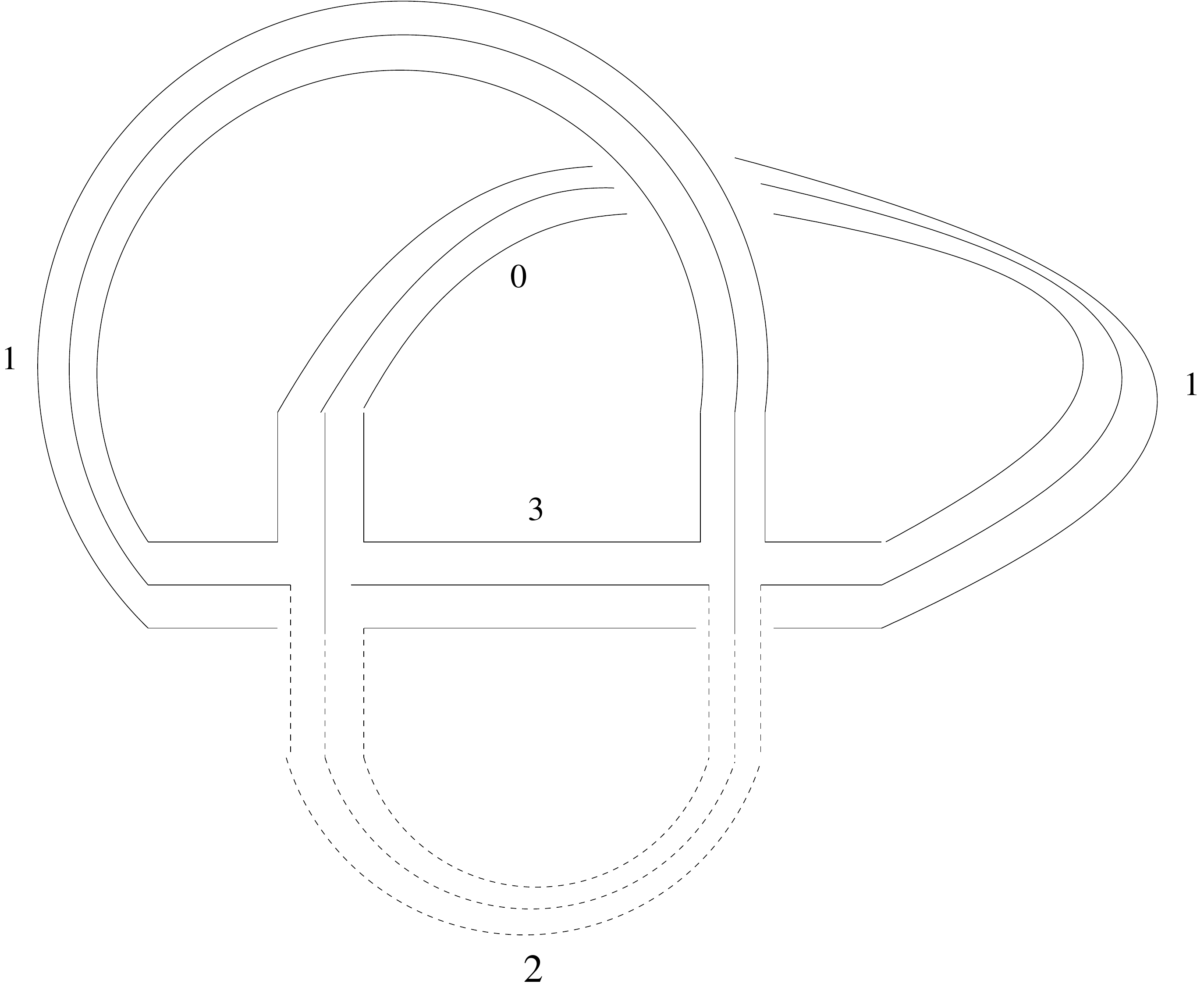}
\caption{The graph $G/e_2$ obtained from the graph $G$ of Fig. \ref{nc} by contraction of the edge $e_2$. This edge has thus become passive, while the edges $e_0, e_1, e_3$ are active edges, as described in section \ref{GFT}.}
\label{Gcontractat}
\end{center}
\end{figure}

Let us now investigate the graph $G/e_2$ which is represented in Fig. \ref{Gcontractat}. 

This graph has three active edges ($e_0, e_1$ and $e_3$) and one passive one ($e_2$, which has just been contracted, as shown in the previous section). The subsets of edges to be investigated are 
\begin{enumerate}
\item  $\{e_0,e_1,e_3\}\cup\{ e_2\}$
\item  $\{e_0,e_1\}\cup\{ e_2\}$,  $\{e_0,e_3\}\cup\{ e_2\}$,  $\{e_1,e_3\}\cup\{ e_2\}$
\item  $\{e_0\}\cup\{ e_2\}$, $\{e_1\}\cup\{ e_2\}$, $\{e_3\}\cup\{ e_2\}$ and
\item $\emptyset\cup\{ e_2\}$.
\end{enumerate}
The first set in this enumeration leads to the contribution $y^3z^2t^2$. The rest of them have no more bubbles and the computation of their contribution follows as for the Bollob\'as-Riordan case. One thus has:
\beqa
\label{del}
{\mathcal T}_{G/e_2}(x,y,z,t)=1+3y+2y^2+y^2z^2+y^3z^2t^2.
\eeqa
Assembling together the two parts \eqref{con} and \eqref{del}, one obtains \eqref{tot}, thus verifying Theorem \ref{teorema}.

The relevance of the graph chosen here comes mainly from the fact that one has a non-trivial dependence in the new $t$ variable, as already mentioned above. The mechanism that we have seen for the emergence of the term(s) containing $t$ does not present any peculiarities for the case of the graph of Fig. \ref{nc}. If we choose to implement the deletion/contraction procedure on the edge $e_0$ (instead of $e_2$) the situation is the analogous. As before, the contribution of the graph $G-e_0$ does not depend on the variable $t$:
\beqa
T_{G-e_0}=y^2+3y+(x-1)+3.
\eeqa
More attention has to be payed for the contracted graph $G/e_0$. The subsets of edges to be investigated are 
$\{e_2,e_1,e_3\}\cup\{ e_0\}$;
$\{e_2,e_1\}\cup\{ e_0\}$,  $\{e_2,e_3\}\cup\{ e_0\}$,  $\{e_1,e_3\}\cup\{ e_0\}$; 
$\{e_2\}\cup\{ e_0\}$, $\{e_1\}\cup\{ e_0\}$, $\{e_3\}\cup\{ e_0\}$ and
 $\emptyset\cup\{ e_0\}$.
As described above, when applying the Definition \ref{def}, the first term in the list above leads to the contribution $y^3z^2t^2$. One has
\beqa
T_{G/e_0}=1+y^3z^2t^2+2y^2z^2+y^2+3y,
\eeqa
thus verifying Theorem \eqref{teorema}.

\section{Conclusions and perspectives}
\label{concluzii}
\renewcommand{\theequation}{\thesection.\arabic{equation}}
\setcounter{equation}{0}

We have thus proposed in this paper a polynomial characterizing tensor graphs; this polynomial makes the difference, in a natural way, between graphs which correspond  to manifold and graphs which do not. Moreover, we have proposed a multivariate and a hypervariate version of this polynomial and we have proved that they satisfy the fundamental deletion/contraction relation.

Since this proposal is made for three-dimensional graphs, an immediate perspective is the generalization of this polynomial to the four-dimensional case. It would then be interesting to understand if there exist any relation between these polynomials and the Krushkal-Renardy one \cite{surpriza}, where, as already mentioned above, a polynomial was introduced for triangulations (without taking into account the tensor graph point of view).

Nevertheless, in order to understand which of these polynomials present interest for theoretical physics, one needs to have some parametric representation for GFT models. An interesting result for combinatorial physics could then be if  this parametric representation can be obtained as some limit of the respective tensor graph polynomial (just as, in commutative scalar QFTs, the parametric representation is obtained as some limit of the Tutte polynomial while, in translation-invariant scalar QFTs on the noncommutative Moyal space, the polynomials of the parametric representation are obtained as some limit of the Bollob\'as-Riordan polynomial \cite{io-BR}).

Let us develop more on this point and its interest for mathematical physics. Considering the multivariate Tutte polynomial, one can obtain under some limit the Symanzik polynomials of the parametric representation of the Feynman amplitudes of the commutative $\Phi^4$ model (for details on the appearance of these polynomials in QFTs, the interested reader can turn to classical textbooks like \cite{iz} or \cite{bookriv}). This relation is possible because the Symanzik polynomials (just like the Tutte polynomial) satisfy the deletion/contraction property (see  \cite{io-BR}, \cite{k} or  \cite{b} for various proofs). 
Thus, is one wants to look for a polynomial which characterizes graph's topology and which is further related to the Symanzik polynomials, the Tutte polynomial is the type of polynomial to look for.

In the case of $\Phi^4$ models on the noncommutative Moyal space, the situation is more complex. Thus, when considering the renormalizable translation-invariant model introduced in \cite{GMRT}, its associated Symanzik polynomials also obey the deletion/contraction relation. This is the case for the multivariate Bollob\'as-Riordan polynomial and one can then prove that the two types of polynomials are related to each others (the Symanzik polynomials of the noncommutative model \cite{GMRT} are obtained as some limit of the  multivariate Bollob\'as-Riordan polynomial).

Moreover, some other types of $\Phi^4$ models are known to be perturbatively renormalizable on the noncommutative Moyal space. This is the case for the Grosse-Wulkenhaar model  \cite{GW} (which, nevertheless is not translation-invariant). In this case, the Symanzik polynomials of its parametric representation do not obey the deletion/contraction relation, but some more involved relations  \cite{KRVT}. It is thus clear that the Symanzik polynomials of this noncommutative QFT model cannot be related to the  Bollob\'as-Riordan polynomial!

Generalizing even further to GFT tensor models, let us advocate that this contraction/relation property can play an important r\^ole for mathematical physics (just like it does for commutative or noncommutative scalar QFTs).

Thus, for the analog of the Symanzik polynomials for some GFT model (some first steps were done in \cite{gft-regain-2} for some simplified version of the four-dimensional topological Ooguri model \cite{ooguri}), it would be interesting to investigate weather or not the deletion/contraction property is satisfied. In this case one knows that, in order to relate these Symanzik polynomials to some polynomial characterizing tensor graphs, polynomials like \eqref{pol} are the type of polynomials to look for.

\medskip

For the sake of completeness, on a more speculative level, let us also recall that the Tutte polynomial is related, under some limit to the Potts model partition function \cite{FK}. It would thus be interesting to explore if some of the generalizations of this Tutte polynomial could also be related to models known in statistical physics.

\section*{Acknowledgments}
 The author acknowledges the grant PN 09 37 01 02 and the CNCSIS grants ``Tinere echipe'' 77/04.08.2010 and  ``Idei'' 454/2009, ID-44. Thomas Krajewski and Vincent Rivasseau are acknowledged for discussions. 
Furthermore, Daniele Oriti is also acknowledged for permission to use Figures \ref{fig:2d} and \ref{fig:t}.

\end{document}